\documentclass{article}
\usepackage[utf8]{inputenc}
\usepackage{amsmath,amstext,amssymb,graphicx,graphics,color}
\usepackage{mathtools}
\usepackage{amsfonts}
\usepackage{theorem}
\usepackage{pgfplots}
\usepackage{cite}
\usepackage{tikz}
\usepackage{hyperref}
\usepackage{bbm}
\usepackage{fancybox}
\usepackage{cleveref}
\usepackage[section]{placeins}
\usepackage{subcaption}

\DeclarePairedDelimiter\floor{\lfloor}{\rfloor}

\theorembodyfont{\normalfont\slshape} 

\newtheorem{theorem}{Theorem}
\newtheorem{proposition}{Proposition}[section]

\newtheorem{lemma}[proposition]{Lemma}
\theoremstyle{break} 

{{\par\noindent \bf Proof. \nobreak}}%
{\nobreak \removelastskip \nobreak \hfill $\Box$ \medbreak}
{{\par\noindent \bf Proof \nobreak}}%
{\nobreak \removelastskip \nobreak \hfill $\Box$ \medbreak}
{{\par\noindent \bf Proof lemma. \nobreak}}%
{\nobreak \removelastskip \nobreak \bf End proof lemma. \medbreak}
\newenvironment{remark}{\par \medskip \noindent {\bf Remark. }\nobreak}{\par \medskip}
\usepackage{fancyhdr}
 \fancypagestyle{plain}{
  \fancyhead{}
  
}

\def\paragraph#1{{\bf #1\ }}
\usepackage{amsmath}
  {
      \theoremstyle{plain}
      
  }
\newcommand{\expo}{\mathrm{e}}

\newcommand{\dd}{\mathrm{d}}

\def\Proof{\noindent{\bf Proof}\quad}
\def\qed{\hfill$\square$\smallskip}

\title{On the equivalence between Fourier-based and Wasserstein distances for probability measures on $\mathbb N$}
\author{Fei Cao\footnotemark[1] \and  Xiaoqian Gong\footnotemark[2]}
\date{\today}

\begin{document}

\maketitle

\footnotetext[1]{University of Massachusetts Amherst - Department of Mathematics and Statistics, 710 N Pleasant St, Amherst, MA 01003, USA}
\footnotetext[2]{Amherst College - Department of Mathematics and Statistics, 220 S Pleasant St, Amherst, MA 01002, USA}

\tableofcontents

\begin{abstract}
In this manuscript we investigate the equivalence of Fourier-based metrics on discrete state spaces with the well-known Wasserstein distances. While the use of Fourier-based metrics in continuous state spaces is ubiquitous since its introduction by Giuseppe Toscani and his colleagues \cite{carrillo_contractive_2007,gabetta_metrics_1995,goudon_fourier_2002} in the study of kinetic-type partial differential equations, the introduction of its discrete analog is recent \cite{bassetti_mean_2015} and seems to be far less studied. In this work, various relations between Fourier-based metrics and Wasserstein distances are shown to hold when the state space is the set of non-negative integers $\mathbb N$. Lastly, we also describe potential applications of such equivalence of metrics in models from econophysics which motivate the present work.
\end{abstract}


\section{Introduction}
\label{sec:sec1}
\setcounter{equation}{0}

Fourier-based metrics, sometimes also known as Toscani distances \cite{rezakhanlou_entropy_2007}, are introduced in a series of works by Giuseppe Toscani and his colleagues \cite{carrillo_contractive_2007,gabetta_metrics_1995,goudon_fourier_2002} for the study of the problem of convergence to equilibrium for the spatially homogenous Boltzmann equation originated from statistical physics. These Fourier-based metrics have also been successfully applied to areas closely related to traditional statistical physics, such as econophysics and sociophysics \cite{cao_binomial_2022,cao_derivation_2021,cao_entropy_2021,cao_explicit_2021,cao_uncovering_2022,during_boltzmann_2008,matthes_steady_2008,naldi_mathematical_2010}. One of the main advantage of using the Fourier-based metrics lies in their interplay with the well-known Wasserstein metrics frequently encountered in the optimal transport literature \cite{villani_optimal_2009}. Indeed, under mild conditions, the equivalence of Fourier-based metrics and Wasserstein metrics (on Euclidean spaces $\mathbb{R}^d$) can be established using elementary analytical and probabilistic tools \cite{carrillo_contractive_2007,gabetta_metrics_1995}. As many partial differential equations (PDEs) such as Boltzmann-type and Fokker-Planck type PDEs have simpler structures under the Fourier coordinates \cite{torregrossa_wealth_2018,toscani_entropy_1999}, it is often quite natural and convenient to transform a given PDE into its equivalent form under Fourier transform and investigate the large time behavior of the transformed PDE by investigating the evolution of Fourier-based metrics with respect to time.

Although the introduction and applications of the Fourier-based metrics when the underlying state space $\mathcal{S}$ is Euclidean (say $\mathcal{S} = \mathbb{R}^d$ with $d\geq 1$) have proven to be very successful and fruitful over the last few decades, its corresponding discrete counterpart, first introduced in \cite{bassetti_mean_2015}, where the state space $\mathcal{S} = \mathbb N$ is countable seems to be far less explored. Indeed, it is still not clear, to the best of our knowledge, whether the Fourier-based metrics between probability measures on $\mathbb N$ will be equivalent to various Wasserstein metrics defined for probability measures supported on $\mathbb N$. Therefore, our primary motivation is to fill in this missing gap and establish such equivalence relation between the aforementioned metrics. It is worth emphasizing that discrete analogues of the classical Fourier-based metrics have been applied in the study of convergence to equilibrium for certain class of infinite system of nonlinear ordinary differential equations (ODEs) \cite{bassetti_mean_2015}, and recently periodic versions of the Fourier-based metrics have seen its application in imaging problems \cite{auricchio_equivalence_2020}.

\subsection{Review of Wasserstein and Fourier-based metrics}

We intend to give a quick and self-contained review of the Wasserstein and Fourier-based metrics, and we refer the interested readers to \cite{carrillo_contractive_2007} for a comprehensive treatments and discussions. We will also restrict ourself to the one-dimensional setting for the ease of presentation, whence $\mathcal{S} = \mathbb{R}$. The Wasserstein distance with exponent $p \geq 1$ for two probability densities $f \in \mathcal{P}(\mathbb R)$ and $g \in \mathcal{P}(\mathbb R)$ (with finite $p$-th order moment) is defined by
\begin{equation}\label{def:Wasserstein}
W_p(f,g) = \left(\inf \int_{\mathbb R \times \mathbb R} |x-y|^2\,\dd \pi(x,y) \right)^{\frac{1}{p}},
\end{equation}
where the infimum runs over all joint densities $\pi$ whose marginal densities coincide with $f$ and $g$, respectively. An alternative formulation for $W_1(f,g)$ involving the cumulative distribution functions turns out to be very useful as well, and it reads as
\begin{equation}\label{def:Wasserstein1}
W_1(f,g) = \int_{\mathbb R} |F(x) - G(x)| \,\dd x,
\end{equation}
where $F$ and $G$ correspond to the cumulative distribution function of $f$ and $g$, respectively. Now we recall the definition of the Fourier-based distance (also known as Toscani's distance) of order $s \geq 1$, given by
\begin{equation}\label{def:Toscani}
d_s(f,g) = \sup\limits_{\xi \in \mathbb{R}\setminus \{0\}} \frac{|\hat{f}(\xi)-\hat{g}(\xi)|}{|\xi|^s},
\end{equation}
where $\hat{f}(\xi) \coloneqq \int_{\mathbb R} f(x)\,\expo^{-i\,x\,\xi}\,\dd x$ denotes the Fourier transform of $f$. It is shown in \cite{carrillo_contractive_2007} that $d_s(f,g) < \infty$ as long as $f$ and $g$ share the same moments up to order $\floor{s}$, where $\floor{s}$ denotes the integer part of $s >0$. Now we recall the following result which ensures the equivalence between $W_2$ and $d_2$ in one space dimension under mild conditions \cite{carrillo_contractive_2007}.

\begin{lemma}\label{lem:classical}
Assume that $f,g \in \mathcal{P}(\mathbb R)$ are probability densities on $\mathbb R$ with equal mean value such that
\[M_{2+\alpha} \coloneqq \max\left\{\int_{\mathbb R} |x|^{2+\alpha}\,f(x)\,\dd x, \int_{\mathbb R} |x|^{2+\alpha}\,g(x)\,\dd x\right\} < \infty\]
for some $\alpha > 0$. Let $C_\alpha = 2^{(2+\alpha)\slash (1+\alpha)}\,\left(\alpha^{1\slash (1+\alpha)}+\alpha^{-\alpha \slash (1+\alpha)}\right)\left(\frac{18}{\pi}\right)^{\alpha\slash \left(3\,(1+\alpha)\right)}$, where \[M_2 \coloneqq \max\left\{\int_{\mathbb R} |x|^2\,f(x)\,\dd x, \int_{\mathbb R} |x|^2\,g(x)\,\dd x\right\} < \infty,\] then we have
\begin{equation}\label{eq:d2_to_W2}
W^2_2(f,g) \leq C_\alpha\,\left[d_2(f,g)\right]^{\alpha\slash \left(6\,(1+\alpha)\right)}\,M^{\alpha\slash \left(3\,(1+\alpha)\right)}_2\,M^{\alpha\slash (1+\alpha)}_{2+\alpha}
\end{equation}
and
\begin{equation}\label{eq:W2_to_d2}
d_2(f,g) \leq \frac 12\,W^2_2(f,g) + \min\left\{\int_{\mathbb R} |x|^2\,f(x)\,\dd x, \int_{\mathbb R} |x|^2\,g(x)\,\dd x\right\}^{\frac 12}\,W_2(f,g)
\end{equation}
\end{lemma}

\Proof We sketch the proof of Lemma \ref{lem:classical} for the sake of completeness, following the argument provided in \cite{carrillo_contractive_2007}. We start with the proof of \eqref{eq:d2_to_W2}. If $X$ is a random variable whose law is $f$, then for each $R > 0$ we have $\mathbb{P}\left[X\geq R\right] \leq \mathbb{E}\left[X^{2+\alpha}\right] \slash R^{2+\alpha}$ thanks to Markov's inequality. This implies \[\lim_{R \to \infty} R^2\,\left(F(-R)+1-F(R)\right) = \lim_{R \to \infty} R^2\,\mathbb{P}\left[|X| \geq R\right] = 0 .\] Therefore, we can integrate by parts the inequalities
\[\int_{-\infty}^{-R} |x|\,|F(x)-G(x)|\,\dd x \leq -\int_{-\infty}^{-R} x\,(F(x)+G(x))\,\dd x \] and
\[\int_R^\infty |x|\,|F(x)-G(x)|\,\dd x \leq \int_R^\infty x\,\left(1-F(x)+1-G(x)\right)\,\dd x\] to obtain
\begin{equation}\label{eq:R1}
\int_{|x|\geq R} |x|\,|F(x)-G(x)|\,\dd x \leq \int_{|x|\geq R} \frac{|x|^2}{2}\,\left(f(x)+g(x)\right)\,\dd x.
\end{equation}
Thus, we deduce that
\begin{align*}
\int_{\mathbb R} |F(x) - G(x)| \,\dd x &\leq \int_{-R}^R |F(x) - G(x)| \,\dd x + \frac{1}{R}\,\int_{|x|\geq R} |x|\,|F(x)-G(x)|\,\dd x \\
&\leq (2\,R)^{\frac 12}\,\left(\int_{\mathbb R} |F(x) - G(x)|^2 \,\dd x\right)^{\frac 12} + \frac{M_2}{R}.
\end{align*}
Optimizing over $R$ gives rise to
\begin{equation}\label{eq:R2}
\int_{\mathbb R} |F(x)-G(x)|\,\dd x \leq \left(2^{-\frac 13}+2^{\frac 23}\right)\,M^{\frac 13}_2\,\left(\int_{\mathbb R} |F(x) - G(x)|^2 \,\dd x\right)^{\frac 13}.
\end{equation}
By the Parseval–Plancherel identity, we also have that
\begin{equation}\label{eq:PP}
\int_{\mathbb R} |F(x) - G(x)|^2 \,\dd x = \frac{1}{2\,\pi}\,\int_{\mathbb R} |\hat{F}(\xi) - \hat{G}(\xi)|^2\,\dd \xi.
\end{equation}
As $\hat{F}(\xi) - \hat{G}(\xi) = \frac{\hat{f}(\xi) - \hat{g}(\xi)}{i\,\xi}$, we estimate
\begin{align*}
\int_{\mathbb R} |F(x) - G(x)|^2 \,\dd x &= \frac{1}{2\,\pi}\,\int_{\mathbb R} \frac{|\hat{f}(\xi) - \hat{g}(\xi)|^2}{\xi^2}\,\dd \xi \\
&\leq \frac{1}{2\,\pi}\,\left(\int_{-R}^R \frac{|\hat{f}(\xi) - \hat{g}(\xi)|^2}{\xi^2}\,\dd \xi + \int_{|\xi|\geq R} \frac{1}{\xi^2}\,\dd \xi \right) \\
&\leq \frac{1}{2\,\pi}\,\left(\int_{-R}^R \xi^2\,d^2_2(f,g)\,\dd \xi + \frac{2}{R} \right) = \frac{1}{\pi}\,\left(\frac{R^3}{3}\,d^2_2(f,g) + \frac{1}{R}\right).
\end{align*}
Optimizing over the choice of $R >0$ once again yields the bound
\begin{equation}\label{eq:R3}
\int_{\mathbb R} |F(x) - G(x)|^2 \,\dd x \leq \frac{4}{3\,\pi}\,\sqrt{d_2(f,g)}.
\end{equation}
Assembling the bounds \eqref{eq:R2} and \eqref{eq:R3} leads us to
\begin{equation}\label{eq:d2_to_W1}
W_1(f,g) \leq \left(\frac{18\,M_2}{\pi}\right)^{\frac 13}\,\left[d_2(f,g)\right]^{\frac 16}.
\end{equation}
By a standard comparison between Wasserstein distances with different exponents, we can also deduce that
\begin{equation}\label{eq:W1_to_W2}
W^2_2(f,g) \leq 2^{(2+\alpha)\slash (1+\alpha)}\,\left(\alpha^{1\slash (1+\alpha)}+\alpha^{-\alpha \slash (1+\alpha)}\right)\,M^{1\slash (1+\alpha)}_{2+\alpha}\,\left[W_1(f,g)\right]^{\alpha \slash (1+\alpha)}.
\end{equation}
As a result, the bound \eqref{eq:d2_to_W2} follows immediately from \eqref{eq:d2_to_W1} and \eqref{eq:W1_to_W2}. To justify the other bound \eqref{eq:W2_to_d2}, let $\pi^*(f,g)$ be the optimal transference plan between $f$ and $g$ for the Euclidean Wasserstein distance $W_2$. As $f$ and $g$ have the same mean value, we have $\int_{\mathbb R \times \mathbb R} (x-y)\,\dd \pi^*(f,g) = 0$, whence
\[\hat{f}(\xi) - \hat{g}(\xi) = \int_{\mathbb R \times \mathbb R} \left(\expo^{-i\,x\,\xi} - \expo^{-i\,y\,\xi} + i\,\xi\,(x-y)\right)\,\dd \pi^*(f,g).\] The integrand can be estimated as
\begin{align*}
|\expo^{-i\,x\,\xi} - \expo^{-i\,y\,\xi} + i\,\xi\,(x-y)| &\leq \left|\expo^{-i\,y\,\xi}\,\left(\expo^{-i\,(x-y)\,\xi}-1+i\,\xi\,(x-y)\right)\right| + |(\expo^{-i\,y\,\xi}-1)\,i\,\xi\,(x-y)| \\
&\leq \frac 12\,|\xi\,(x-y)|^2 + |y|\,|\xi|\,|\xi\,(x-y)| \\
&= \frac 12\,|\xi|^2\,|x-y|^2 + |\xi|^2\,|y|\,|x-y|.
\end{align*}
Integrating the previous inequality yields
\begin{align*}
\frac{|\hat{f}(\xi)-\hat{g}(\xi)|}{|\xi|^2} &\leq \int_{\mathbb R \times \mathbb R} \left(\frac 12\,|x-y|^2 + |y|\,|x-y|\right)\,\dd \pi^*(f,g)\\
&\leq \frac 12\,W^2_2(f,g) + \left(\int_{\mathbb R} |x|^2\,g(x)\,\dd x\right)^{\frac 12}\,W_2(f,g)
\end{align*}
for all $\xi \in \mathbb{R} \setminus \{0\}$. Finally, the advertised bound \eqref{eq:W2_to_d2} follows from an obvious symmetrization of the previous inequality. \qed


\subsection{Toscani and Wasserstein metrics for probabilities on $\mathbb N$}

As has been mentioned in the beginning of the manuscript, the natural analog of the classical Fourier-based distance between probability measures over the set of non-negative integers (i.e., $\mathcal{S} = \mathbb{N}$) is introduced in a recent work \cite{bassetti_mean_2015}, and its $s$-th order distance (with $s \geq 1$) is defined via
\begin{equation}\label{def:Toscani_pdf}
D_s(f,g) = \sup\limits_{z \in (0,1)} \frac{|\hat{f}(z)-\hat{g}(z)|}{|1-z|^s},
\end{equation}
where $f = (f_0,f_1,\ldots)$, $g = (g_0,g_1,\ldots)$ are elements of $\mathcal{P}(\mathbb N)$ and $\hat{f}(z)$ represents the probability generating function associated with the law $f$, i.e.,
\begin{equation}\label{def:pdf}
\hat{f}(z) \coloneqq \sum\limits_{n \geq 0} z^n\,f_n, \quad z \in [0,1].
\end{equation}
We remark here that the function $\hat{f}(z)$ is termed as the probability generating function due to the following easy-to-verify property:
\begin{equation*}
\mathbb{P}\left[X = n\right] = f_n = \frac{\hat{f}^{(n)}(z)}{n!}\Big|_{z=0},
\end{equation*}
where $X$ is a random variable distributed according to $f$. As was already observed in \cite{bassetti_mean_2015}, convergence with respect to the distance $D_s$ implies convergence of probability generating functions, which in turn yields the pointwise convergence of probability distributions on $\mathbb{N}$. We also recall that, for a pair of probability distributions $f \in \mathcal{P}(\mathbb N)$ and $g \in \mathcal{P}(\mathbb N)$, the Wasserstein distance of order $p \geq 1$ between $f$ and $g$ is defined via
\begin{equation}\label{def:Wasserstein_discrete}
W_p(f,g) = \left(\inf\limits_{\pi} \sum_{i\geq 0}\sum_{j\geq 0} \pi_{i,j}\,|i-j|^p ~\bigg|~ \sum\limits_{i\geq 0} \pi_{i,j} = g_j,~ \sum\limits_{j\geq 0} \pi_{i,j} = f_i,~\pi_{i,j} \geq 0 \right)^{\frac{1}{p}}.
\end{equation}
To the best of our knowledge, the relation between the Toscani distance \eqref{def:Toscani_pdf} and the Wasserstein distance \eqref{def:Wasserstein_discrete} (when the underlying state space is the set of non-negative integers $\mathbb N$) still remains to be uncovered. This manuscript therefore aims at proving the equivalence between the aforementioned metrics, and in particular, we will show that

\begin{theorem}\label{thm:main}
Assume that $f,g\in \mathcal{P}(\mathbb N)$, we have
\begin{equation}\label{part1}
D_1(f,g) \leq W_1(f,g).
\end{equation}
If $f \in \mathcal{P}(\mathbb N)$ and $g \in \mathcal{P}(\mathbb N)$ have the same mean value, then
\begin{equation}\label{part2}
D_2(f,g) \leq \frac 12\,W^2_2(f,g) + \min\left\{\sum\limits_{n\geq 0} n^2\,f_n, \sum\limits_{n\geq 0} n^2\,g_n\right\}^{\frac 12}\,W_2(f,g).
\end{equation}
Lastly, assume that $f,g\in \mathcal{P}(\mathbb N)$ with equal mean value such that $\mathrm{supp}(f) \cup \mathrm{supp}(g) \subset \{0,1\ldots,N\}$ for some finite $N \in \mathbb N_+$. Assume also that
\[m_{2+\alpha} \coloneqq \max\left\{\sum\limits_{n\geq 0} n^{2+\alpha}\,f_n,\sum\limits_{n\geq 0} n^{2+\alpha}\,g_n\right\} < \infty\]
for some $\alpha > 0$. Then we have
\begin{equation}\label{part3}
W^2_2(f,g) \leq C\,\left[D_2(f,g)\right]^{\alpha \slash (1+\alpha)}.
\end{equation}
for some $C = C(N,\alpha,m_{2+\alpha}) > 0$ depending only on $N$, $\alpha$ and $m_{2+\alpha}$.
\end{theorem}

The rest of the manuscript will be fully devoted to the proof of Theorem \ref{thm:main} and is organized as follows: we establish the first part \eqref{part1} of Theorem \ref{thm:main} in section \ref{subsec:1}, allowing us to bound $D_1(f,g)$ by $W_1(f,g)$. Sections \ref{subsec:2} and \ref{subsec:3} are devoted to the second part \eqref{part2} and the third part \eqref{part3} of Theorem \ref{thm:main}, respectively. We demonstrate in section \ref{sec:sec3} an important application that motivates the study of equivalence relations between Fourier-based metrics and Wasserstein distances in the discrete setting. Finally, we draw a brief conclusion in section \ref{sec:sec4} and also provide several potential directions for future research activities.

\section{Equivalence between Toscani and Wasserstein distances}
\label{sec:sec2}
\setcounter{equation}{0}

For notational simplicity, we denote $i \wedge j \coloneqq \min\{i,j\}$ and $i \vee j \coloneqq \max\{i,j\}$ for each pair of non-negative integers $(i,j) \in \mathbb{N}\times \mathbb{N}$.

\subsection{From $W_1$ to $D_1$}\label{subsec:1}

We establish the following bound which allows us to control the Toscani distance $D_1$ by the Wasserstein distance $W_1$. We emphasize that the same relation holds when the relevant distances are defined on Euclidean spaces \cite{carrillo_contractive_2007,during_boltzmann_2008}.

\begin{lemma}[From $W_1$ to $D_1$]\label{lem:2.1}
Assume that $f,g\in \mathcal{P}(\mathbb N)$, we have
\begin{equation}\label{eq:W1_to_D1}
D_1(f,g) \leq W_1(f,g).
\end{equation}
\end{lemma}

\Proof Let $\pi^*(f,g)$ be the optimal transference plan between $f$ and $g$ for the Wasserstein distance $W_1$. We observe that
\begin{equation}\label{eq:OBS1}
\frac{|z^i - z^j|}{1-z} = \left|\sum\limits_{n = i\wedge j}^{i\vee j - 1} z^n \right| \leq |i-j|
\end{equation}
for each $(i,j) \in \mathbb{N}\times \mathbb{N}$ and for all $z \in (0,1)$. Therefore, we deduce that
\begin{align*}
D_1(f,g) &= \sup\limits_{z \in (0,1)} \frac{1}{1-z}\,\left|\sum_{i\geq 0} z^i\,f_i - \sum_{j \geq 0} z^j\,g_j\right| \\
&= \sup\limits_{z \in (0,1)} \frac{1}{1-z}\,\left|\sum_{i\geq 0}\sum_{j\geq 0} (z^i-z^j)\,\pi^*_{i,j}\right| \\
&\leq \sup\limits_{z \in (0,1)} \sum_{i\geq 0}\sum_{j\geq 0} \frac{|z^i - z^j|}{1-z}\,\pi^*_{i,j} \\
&\leq \sum_{i\geq 0}\sum_{j\geq 0} |i-j|\,\pi^*_{i,j} = W_1(f,g),
\end{align*}
whence the proof is completed.

\begin{remark}
An alternative route to justify the content of Lemma \ref{lem:2.1} relies on the following expression for the Wasserstein distance $W_1$:
\begin{equation}\label{eq:dual_W1}
W_1(f,g) = \sum\limits_{n\geq 0} \left|\sum\limits_{k=0}^n f_k - \sum\limits_{k=0}^n g_k\right|. 
\end{equation}
Indeed, for each $z \in (0,1)$ we have
\begin{equation}\label{eq:chain_of_inequalities}
\begin{aligned}
\left|\frac{\hat{f}(z)-\hat{g}(z)}{1-z}\right| &= \left|\sum_{k\geq 0} z^k\,(f_k-g_k) \cdot \sum_{\ell \geq 0} z^\ell \right| \\
&= \left|\sum_{k\geq 0}\sum_{\ell \geq 0} z^{k+\ell} \,(f_k-g_k) \right| \\
&= \left|\sum_{n\geq 0} z^n\,\sum_{k=0}^n (f_k-g_k)\right| \leq W_1(f,g).
\end{aligned}
\end{equation}
Taking the supremum of the left hand side of \eqref{eq:chain_of_inequalities} over $z \in (0,1)$ yields the claimed bound \eqref{eq:W1_to_D1}.
\end{remark}

\subsection{From $W_2$ to $D_2$}\label{subsec:2}

Next, we justify the domination of the Toscani distance $D_2$ by the Wasserstein distance $W_2$. Once again, the relation to be shown below shares the same spirit as the case when the relevant distances are considered on Euclidean spaces \cite{carrillo_contractive_2007,during_boltzmann_2008}.

\begin{lemma}[From $W_2$ to $D_2$]\label{lem:2.2}
Assume that $f,g\in \mathcal{P}(\mathbb N)$ with equal mean value, we have
\begin{equation}\label{eq:W2_to_D2}
D_2(f,g) \leq \frac 12\,W^2_2(f,g) + \min\left\{\sum\limits_{n\geq 0} n^2\,f_n, \sum\limits_{n\geq 0} n^2\,g_n\right\}^{\frac 12}\,W_2(f,g).
\end{equation}
\end{lemma}

\Proof Let $\pi^*(f,g)$ be the optimal transference plan between $f$ and $g$ for the Wasserstein distance $W_2$. As $f$ and $g$ share the same mean value, we have
\begin{equation*}
\sum\limits_{i\geq 0}\sum\limits_{j\geq 0} (i-j)\,\pi^*_{i,j} = 0.
\end{equation*}
We also observe, by a simple induction argument, that
\begin{equation}\label{eq:OBS2}
\left|z^m - 1 + m\,(1-z)\right| \leq \frac{m^2\,(1-z)^2}{2} \quad \textrm{and} \quad 1 - z^m \leq m\,(1-z)
\end{equation}
for all $m \in \mathbb N$ and all $z \in (0,1)$. Now for each $z \in (0,1)$ we estimate $|\hat{f}(z) - \hat{g}(z)|$ as
\begin{align*}
\left|\hat{f}(z)-\hat{g}(z)\right|&=\left|\sum\limits_{i\geq 0}\sum\limits_{j\geq 0} \left(z^i - z^j + (i-j)\,(1-z)\right)\,\pi^*_{i,j}\right| \\
&\leq \sum\limits_{i\geq 0}\sum\limits_{j\geq 0} \left|\left(z^i - z^j + (i-j)\,(1-z)\right)\right|\,\pi^*_{i,j} \\
&\leq \sum\limits_{i\geq 0}\sum\limits_{j\geq 0} z^{i\wedge j}\,\left|z^{i\vee j - i\wedge j} - 1 + (i\vee j - i\wedge j)\,(1-z)\right|\,\pi^*_{i,j} \\
&\qquad + \sum\limits_{i\geq 0}\sum\limits_{j\geq 0} \left(1-z^{i\wedge j}\right)\,(i\vee j - i\wedge j)\,(1-z)\,\pi^*_{i,j} \\
&\leq \sum\limits_{i\geq 0}\sum\limits_{j\geq 0} \left(\frac{|i-j|^2}{2}\,(1-z)^2 + (i\wedge j)\,|i-j|\,(1-z)^2\right)\,\pi^*_{i,j}.
\end{align*}
Consequently, for each $z \in (0,1)$ we obtain
\[\frac{\left|\hat{f}(z) - \hat{g}(z)\right|}{(1-z)^2} \leq \sum\limits_{i\geq 0}\sum\limits_{j\geq 0} \left(\frac{|i-j|^2}{2} + (i\wedge j)\,|i-j|\right)\,\pi^*_{i,j},\]
from which the advertised bound \eqref{eq:W2_to_D2} follows by Cauchy-Schwarz inequality and the definition of the Wasserstein distance $W_2$. \qed

\subsection{From $D_2$ to $W_2$}\label{subsec:3}

Finally, we take on the hardest part of the proof of Theorem \ref{thm:main}, which requires a control of the Wasserstein distance $W_2$ by the Toscani distance $D_2$. Indeed, since we no longer have access to a natural analog of the Parseval–Plancherel identity involving the probability generating function, we have to discard the proof strategy outlined in the proof of Lemma \ref{lem:classical} and resort to new tools. We now establish an upper bound of $W_2$ in terms of $D_2$, at least when the support of the probability distributions are restricted to a compact subset of $\mathbb N$.

\begin{lemma}[From $D_2$ to $W_2$ for compactly supported distributions]\label{lem:2.3}
Assume that $f,g\in \mathcal{P}(\mathbb N)$ with equal mean value such that $\mathrm{supp}(f) \cup \mathrm{supp}(g) \subset \{0,1\ldots,N\}$ for some finite $N \in \mathbb N_+$. Assume also that
\[m_{2+\alpha} \coloneqq \max\left\{\sum\limits_{n\geq 0} n^{2+\alpha}\,f_n,\sum\limits_{n\geq 0} n^{2+\alpha}\,g_n\right\} < \infty\]
for some $\alpha > 0$. Then we have
\begin{equation}\label{eq:D2_to_W2}
W^2_2(f,g) \leq C\,\left[D_2(f,g)\right]^{\alpha \slash (1+\alpha)}.
\end{equation}
for some $C = C(N,\alpha,m_{2+\alpha}) > 0$ depending only on $N$, $\alpha$ and $m_{2+\alpha}$.
\end{lemma}

\Proof We let $F$ and $G$ to be the cumulative distribution functions corresponding to $f$ and $g$, respectively. In other words, $F_n = \sum_{k=0}^n f_k$ and $G_n = \sum_{k=0}^n g_k$ for $n \in \mathbb N$. By our assumption on the support of $f$ and $g$, $F_n = G_n = 1$ for all $n\geq N$ or equivalently $f_n = g_n = 0$ for all $n\geq N+1$. Therefore, we can write
\[W_1(f,g) = \sum\limits_{n\geq 0} |F_n - G_n| = \sum\limits_{n=0}^N |F_n - G_n|.\] Since $f$ and $g$ have the same mean value, i.e., $\sum_{n=0}^N n\,f_n = \sum_{n=0}^N n\,g_n$, we deduce that $\sum_{n=0}^N (F_n - G_n) = 0$ or equivalently $F_0 - G_0 = -\sum_{n=1}^N (F_n - G_n)$. As a consequence, we arrive at
\begin{equation}\label{eq:b1}
W_1(f,g) = \sum\limits_{n=0}^N |F_n - G_n| \leq 2\,\sum\limits_{n=1}^N |F_n - G_n|.
\end{equation}
To bound $W_1(f,g)$ from above by $D_2(f,g)$ we will make a digression first. We claim that
\begin{equation}\label{eq:ell}
\begin{aligned}
\mathbb R^{N} &\to [0, \infty)\\
\ell \colon \boldsymbol{a} &\mapsto \sup\limits_{z\in (0,1)}\left|\sum_{n=1}^N a_n\sum_{k=0}^{n-1} z^k\right|
\end{aligned}
\end{equation}
defines a norm on $\mathbb{R}^N$. Indeed, it is clear from the definition that $\ell[\boldsymbol{a}] \geq 0$ and $\ell[c\cdot\boldsymbol{a}] = |c|\,\ell[\boldsymbol{a}]$ for all $c \in \mathbb R$. Assume that $\ell[\boldsymbol{a}] = 0$, then for all $z\in (0, 1)$, \[\sum_{n=1}^N a_n \sum_{k=0}^{n-1} z^k = \sum_{k=0}^{N-1} \left(\sum_{n=k+1}^{N} a_n\right)z^k = 0,\] which implies that $\sum_{n=k+1}^N a_n = 0$ for all $0\leq k\leq N-1$. Thus $a_N = 0$ and by induction on $k = 1,\ldots,N-1$, $a_k  = \sum_{n=k}^{N}a_n - \sum_{n=k+1}^N a_n  = 0$. Therefore, $\boldsymbol{a} = 0$. Lastly, we check the triangle inequality. For $\boldsymbol{a}$ and $\boldsymbol{b}$ in $\mathbb R^N$, as
\begin{align*}
\left|\sum_{n=1}^{N} \left(a_n + b_n\right) \sum_{k=0}^{n-1} z^k\right| &= \left|\sum_{n=1}^{N} a_n \sum_{k=0}^{n-1} z^k + \sum_{n=1}^{N} b_n \sum_{k=0}^{n-1} z^k\right| \\
&\leq \left|\sum_{n=1}^{N} a_n \sum_{k=0}^{n-1} z^k\right| + \left|\sum_{n=1}^{N} b_n \sum_{k=0}^{n-1} z^k\right|
\end{align*}
for all $z\in (0, 1)$, $\ell\left[\boldsymbol{a} + \boldsymbol{b}\right] \leq \ell[\boldsymbol{a}] + \ell[\boldsymbol{b}]$ and we conclude that $\ell \colon \mathbb R^{N} \to [0, \infty)$ is a norm on $\mathbb{R}^N$. In particular, as all norms on $\mathbb{R}^N$ are equivalent, there exists some universal constant $C = C(N) > 0$ such that
\begin{equation}\label{eq:equivalence}
\|\boldsymbol{a}\|_{\ell^1} = \sum\limits_{n=1}^N |a_n| \leq C\,\ell[\boldsymbol{a}].
\end{equation}
Now, we notice that \[\hat{f}(z) = \sum\limits_{n=0}^N z^n\,f_n = \sum\limits_{n=0}^N z^n\,(F_n - F_{n-1}) = \sum\limits_{n=0}^N (z^n - z^{n+1})F_n + z^{N+1}\] and similarly $\hat{g}(z) = \sum\limits_{n=0}^N (z^n - z^{n+1})G_n + z^{N+1}$. Thus,
\begin{equation}\label{eq:chain}
\begin{aligned}
\frac{|\hat{f}(z)-\hat{g}(z)|}{(1-z)^2} &= \frac{1}{1-z}\,\left|\sum\limits_{n=0}^N z^n\,(F_n - G_n)\right| \\
&= \frac{1}{1-z}\,\left|\sum\limits_{n=1}^N (z^n-1)\,(F_n - G_n)\right|
&= \left|\sum\limits_{n=1}^N (F_n - G_n)\,\sum\limits_{k=1}^{N-1} z^k\right|.
\end{aligned}
\end{equation}
Applying the relation \eqref{eq:equivalence} with $a_n = F_n - G_n$ and employing the bound \eqref{eq:b1}, we deduce that
\begin{equation}\label{eq:b2}
W_1(f,g) \leq 2\,C\,D_2(f,g).
\end{equation}
Since we can still bound $W_2(f,g)$ by $W_1(f,g)$ as in \eqref{eq:W1_to_W2} in the discrete setting (except that $M_{2+\alpha}$ in \eqref{eq:W1_to_W2} needs to be replaced by $m_{2+\alpha}$), the advertised inequality \eqref{eq:D2_to_W2} follows and the proof of Lemma \ref{lem:2.3} is completed. \qed

\begin{remark}
Unfortunately, our Lemma \ref{lem:2.3} can only hand the case of compactly supported probability distributions as the constant $C = C(N,\alpha,m_{2+\alpha})$ appearing in \eqref{eq:D2_to_W2} might blow up to infinity as $N \to \infty$. Due to the lack of Parseval–Plancherel identity involving probability generating functions, we are unable treat the general case of probability distributions with full support on $\mathbb N$, and we leave it as a remaining open problem.
\end{remark}




\section{Application to econophysics models}
\label{sec:sec3}
\setcounter{equation}{0}

Econophysics is a sub-branch of statistical physics that applies traditional concepts and techniques from classical physics to models of economical systems. We refer the interested readers to the review article \cite{dragulescu_statistical_2000} where many models from econophysics were first introduced, and we also refer to \cite{cao_binomial_2022,cao_derivation_2021,cao_entropy_2021,cao_explicit_2021,cao_uncovering_2022,cao_interacting_2022,cortez_uniform_2022,lanchier_rigorous_2018,lanchier_rigorous_2019,lanchier_distribution_2022} for many other models arising from econophysics literature.

The primary motivation behind this work is motivated by the so-called binomial reshuffling model for money exchange in a closed economic system with $N$ agents, introduced in a very recent work \cite{cao_binomial_2022} which we intend to briefly describe the set-up. The binomial reshuffling dynamics consists in choosing at random time two individuals and to redistribute their money between them according to a binomial distribution. To write this dynamics rigorously, we denote by $X_i(t) \in \mathbb N$ the amount of dollars the agent $i$ has at time $t$ and we suppose that $\sum_{i=1}^N X_i(0) = N\,\mu$ for some fixed $\mu \in \mathbb{N}_+$. At a random time generated by a Poisson clock with rate $N$, two agents (say $i$ and $j$) update their wealth according to the following rule:
\begin{equation}
\label{binomial_reshuffling}
(X_i,X_j)~ \begin{tikzpicture} \draw [->,decorate,decoration={snake,amplitude=.4mm,segment length=2mm,post length=1mm}]
(0,0) -- (.6,0); \end{tikzpicture}~  \left(B\circ (X_i + X_j), X_i + X_j - B\circ (X_i + X_j)\right),
\end{equation}
where $B\circ (X_i + X_j)$ represents a binomial random variable with parameter $(X_i+X_j,\frac 12)$. The agent-based numerical simulation suggests that, as the number of agents and time go to infinity, the limiting distribution of money is well-approximated by a Poisson distribution with parameter $\mu$.

If we denote by ${\bf p}(t)=\left(p_0(t),p_1(t),\ldots,p_n(t),\ldots\right)$ the law of the process $X_1(t)$ as $N \to \infty$, i.e., $p_n(t) = \lim\limits_{N \to \infty} \mathbb P\left(X_1(t) = n\right)$, then a heuristic mean-field type argument as carried out in \cite{cao_binomial_2022} shows that the time evolution of ${\bf p}(t)$ is given by
\begin{equation}\label{eq:ODE}
\frac{\dd}{\dd t} {\bf p}(t) = Q[{\bf p}(t)]
\end{equation}
where
\begin{equation}\label{eq:Q}
Q[{\bf p}]_n = \sum\limits_{k\geq 0}\sum\limits_{\ell \geq 0} \tbinom{k+\ell}{n}\,\tfrac{1}{2^{k+\ell}}\,p_k\,p_\ell\,\mathbbm{1}\{n \leq k+\ell\} - p_n \quad \forall~ n \geq 0,
\end{equation}
with the usual convention that $\tbinom{0}{0}$ is interpreted as $0$. Moreover, the dynamics preserves the mean value $\sum_{n\geq 0} n\,p_n$ and it can be shown that the Poisson distribution ${\bf p}^* = \{p^*_n\}_{n\in \mathbb N}$ defined by
\begin{equation}\label{eq:equil}
p^*_n = \frac{\mu^n\,\expo^{-\mu}}{n!}, \quad n \in \mathbb N
\end{equation}
is the unique equilibrium distribution associated with the mean-field ODE system \eqref{eq:ODE}. The main result in \cite{cao_binomial_2022} is the following quantitative convergence guarantee measured in terms of the $W_2$ Wasserstein distance: there exists some constance $C > 0$ depending only on $W_2({\bf p}(0),{\bf p}^*)$ such that
\begin{equation}\label{eq:Wasserstein_convergence}
W_2({\bf p}(t),{\bf p}^*) \leq \frac{C}{\sqrt{t}},~~\forall t > 0.
\end{equation}
On the other hand, the nonlinear ODE system \eqref{eq:ODE} and \eqref{eq:Q} is a special case of a more general ODE system motivated by biological applications and investigated by Federico Bassetti and Giuseppe Toscani in \cite{bassetti_mean_2015}, where they established a exponential decay result in terms of the $D_2$ distance: there exists some constant $\alpha > 0$ depending only on second moment of the initial datum ${\bf p}(0)$ such that
\begin{equation}\label{eq:D2_convergence}
D_2({\bf p}(t),{\bf p}^*) \leq D_2({\bf p}(0),{\bf p}^*)\,\expo^{-\alpha\,t},~~\forall t > 0.
\end{equation}
Therefore, it is very natural to ask whether one can translate the exponential decay in $D_2$ distance towards an exponential decay in $W_2$ distance (observed numerically in \cite{cao_binomial_2022} for some initial conditions). If one can strengthen the content of Lemma \ref{lem:2.3} to handle the general case of probability distributions with full support on $\mathbb N$, so that a bound in the same spirit as \eqref{eq:D2_to_W2} still remains valid, then our previous question will have a positive and affirmative answer. However, we fail to carry out the desired refinement of Lemma \ref{lem:2.3} in this work.

\section{Conclusion}
\label{sec:sec4}
\setcounter{equation}{0}

In this manuscript, we investigated various relations between the Fourier-based metrics and the well-known Wasserstein distances in a discrete setting. While applications of Fourier-based metrics introduced by Giuseppe Toscani and his colleagues \cite{auricchio_equivalence_2020,carrillo_contractive_2007,gabetta_metrics_1995,goudon_fourier_2002} to the study of the long time behaviour for Boltzmann-type and Fokker-Planck type PDEs have been well-documented in literature, the introduction and study of such metrics where the state space is discrete \cite{bassetti_mean_2015} have not received enough attention in our opinion. Given the fact that Fourier-based metrics have been applied to a variety of models from econophysics and sociophysics \cite{during_boltzmann_2008,during_continuum_2022,garibaldi_statistical_2007,li_bounded_2020,matthes_steady_2008,naldi_mathematical_2010,torregrossa_wealth_2018,toscani_entropy_1999}, we believe that a detailed investigation of the properties of Fourier-based metrics in discrete settings will benefit the community of kinetic theory in general, and the communities of econophysics and sociophysics in particular. In fact, our primary motivation for the current work is motivated by the binomial reshuffling model introduced in a recent work \cite{cao_binomial_2022} and reviewed in section \ref{sec:sec3}, where the equivalence of Fourier-based metrics and Wasserstein distances can lead to exponential decay result under both metrics for the mean-field version of the model.

As we have emphasized earlier, the lack of Parseval-Plancherel type identity involving the probability generating function is a major technical bottleneck which hinders the attempt to bound Wasserstein distances (from above) by Fourier-based metrics, and therefore we leave it as a open problem for the improvement of the last part of Theorem \ref{thm:main} (equivalently Lemma \ref{lem:2.3}) so that case of probability distributions with full support on $\mathbb N$ can be handled.

Lastly, we comment on several disadvantages of using Fourier-based metrics compared to Wasserstein distances, and suggest a few potential directions for future research. First, it seems to us that Fourier-based metrics have no equivalent formulation (besides its very definition), this is in sharp contrast to Wasserstein distances \cite{villani_optimal_2009}, which admit a duality formulation due to Kantorovich and Rubinstein, a probabilistic interpretation based on coupling, and a kinetic energy minimization formulation due to Benamou and Brenier, just to name a few. Thus it would be quite interesting to investigate whether convenient equivalent definitions of Fourier-based metrics exist or not so that these metrics can reach out to other communities in pure and/or applied mathematics. Second, it seems that simulation of Fourier-based metrics is much harder compared to the numerical implementation of Wasserstein metrics, especially given the fact that the computational aspects of Wasserstein metrics are well-documented \cite{santambrogio_optimal_2015,peyre_computational_2019}. Thus, it would also be interesting to develop convenient numerical tools which enable us to simulate the evolution of Fourier-based metrics in the context of evolution PDEs studied in econophysics and sociophysics literature that we mentioned before.

\vspace{0.2in}

\noindent {\bf Acknowledgement} We highly appreciate the help from a user on Math Stack Exchange where a key component used in the proof of Lemma \ref{lem:2.3} is provided \cite{Kroki_2024}.

\end{document}